\documentclass[letterpaper, 10pt, twocolumn]{article}

\usepackage{sectsty}
\sectionfont{\fontsize{12}{15}\selectfont}

\usepackage[margin={1.5cm,1.5cm}]{geometry}
\usepackage{amssymb,amsmath,color} 
\usepackage{graphicx} 
\usepackage{epsfig}
\usepackage{algorithm}
\usepackage[noend]{algpseudocode}

\usepackage{psfrag}

\newcommand{\real}{{\mathbb{R}}}

\newcommand{\until}[1]{\{1,\ldots,#1\}}

\newcommand{\EE}{\mathcal{E}} 
\newcommand{\GG}{\mathcal{G}}

 \newcommand{\subj}{\text{subj. to}}

\newcommand{\prox}{\mathbf{prox}}
\newcommand{\argmin}{\mathop{\rm argmin}}
\newcommand{\argmax}{\mathop{\rm argmax}}

\newcommand{\DDPGfixed}{Distributed Dual Proximal Gradient}
\newcommand{\DDPGgossip}{Asynchronous Distributed Dual Proximal Gradient}

\newcommand{\nbrs}{\mathcal{N}}
\newcommand{\dom}{\mathop{\bf dom}} 

\algdef{SE}[DOWHILE]{Do}{doWhile}{\hskip\algorithmicindent\algorithmicdo}[1]{\hskip\algorithmicindent\algorithmicwhile\ #1}%

\makeatletter
\newcommand{\StatexIndent}[1][3]{%
  \setlength\@tempdima{\algorithmicindent}%
  \Statex\hskip\dimexpr#1\@tempdima\relax}
\makeatother

\algnewcommand{\algorithmicgoto}{\textbf{go to }}%
\algnewcommand{\Goto}[1]{\algorithmicgoto Line~\ref{#1}}%
\algnewcommand{\Label}{\State\unskip}

\renewcommand{\algorithmicwhile}{\hskip\algorithmicindent \textbf{While:}}

\newtheorem{theorem}{Theorem}[section]

 \newtheorem{lemma}[theorem]{Lemma}
\newtheorem{remark}[theorem]{Remark}

\newtheorem{assumption}[theorem]{Assumption}
 
\newcommand\oprocendsymbol{\hbox{$\square$}}
\newcommand\oprocend{\relax\ifmmode\else\unskip\hfill\fi\oprocendsymbol}

\begin{document}

\title{Randomized dual proximal gradient\\ for large-scale distributed optimization}

\author{Ivano Notarnicola and Giuseppe Notarstefano}
\date{}

\maketitle
    
\let\thefootnote\relax\footnotetext{
  Ivano Notarnicola and Giuseppe Notarstefano are with the Department of Engineering,
    Universit\`a del Salento, Via Monteroni, 73100
    Lecce, Italy, \texttt{name.lastname@unisalento.it.} This result is part of a
    project that has received funding from the European Research Council (ERC)
    under the European Union’s Horizon 2020 research and innovation programme
    (grant agreement No 638992 - OPT4SMART).  } 
    
\begin{abstract}
  In this paper we consider distributed optimization problems in which the cost
  function is separable (i.e., a sum of possibly non-smooth functions all
  sharing a common variable) and can be split into a strongly convex term and a
  convex one. The second term is typically used to encode constraints or to
  regularize the solution. We propose an asynchronous, distributed optimization
  algorithm over an undirected topology, based on a proximal gradient update on
  the dual problem. We show that by means of a proper choice of primal
  variables, the dual problem is separable and the dual variables can be stacked
  into separate blocks. This allows us to show that a distributed
  gossip update can be obtained by means of a randomized block-coordinate
  proximal gradient on the dual function.
\end{abstract}

\section{Introduction}
\label{sec:intro}
Several estimation, learning, decision and control problems arising in
cyber-physical networks involve the distributed solution of a constrained
optimization problem. 
Typically, computing processors have only a partial knowledge of the problem
(e.g. a portion of the cost function or a subset of constraints) and need to
cooperate to compute a global solution of the whole problem.
A key challenge to take into account when designing distributed optimization
algorithms in peer-to-peer networks is that the communication is time-varying
and possibly asynchronous, see, e.g., \cite{tsianos2012consensus} for a review.
%

Early references on distributed optimization algorithms involved primal and dual
subgradient methods and Alternating Direction Method of Multipliers (ADMM),
designed for synchronous communication protocols over fixed graphs. More
recently time-varying versions of these algorithmic ideas have been proposed to
cope with more realistic peer-to-peer network scenarios.
%
%
A Newton-Raphson consensus strategy is proposed in
\cite{zanella2012asynchronous} to solve unconstrained, convex optimization
problems under asynchronous, symmetric gossip communications.
In \cite{jakovetic2014convergence} the authors propose accelerated distributed
gradient methods for unconstrained problems over symmetric, time-varying
networks connected on average.
In order to deal with (time-varying) directed graphs, in
\cite{nedic2013distributed} a push-sum algorithm for average consensus is
combined with a primal subgradient method. %
Paper \cite{akbari2014distributed} extends these methods to online distributed
optimization.
In \cite{kia2014} a novel class of continuous-time, gradient-based distributed
algorithms is proposed. 
A distributed (primal) proximal-gradient method is proposed in
\cite{chen2012fast} for separable optimization problems 
which can handle only a common constraint.
To solve constrained convex optimization problems, in \cite{lee2013distributed}
a distributed random projection algorithm is proposed for a balanced
time-varying graph. %

In \cite{wei20131} a novel asynchronous ADMM-based distributed method is
proposed for separable, constrained convex optimization problem. 
%
%
In \cite{necoara2013random} the author proposes (primal) randomized
block-coordinate descent methods for minimizing multi-agent convex optimization
problems with linearly coupled constraints over networks.
A combination of successive approximations and block-coordinate updates is
proposed in \cite{facchinei2015parallel} to solve separable, non-convex
optimization problems in a big-data setting.
Another class of algorithms exploits the exchange of active constraints among
the nodes to solve general constrained convex programs
\cite{notarstefano2007distributed}. The constraint exchange idea has been
combined with dual decomposition and cutting-plane methods to solve robust
convex optimization problems via polyhedral approximations
\cite{burger2014polyhedral}. These algorithms work under asynchronous, directed
and unreliable communication.

It is worth noting that the algorithm in \cite{necoara2013random} uses a
coordinate-descent idea similar to the one we use in this paper, but it works
directly on the primal problem. Similarly, in \cite{chen2012fast} the proximal
operator is used to handle the sparsity constraints directly on the primal
problem, so that local constraints cannot be simultaneously taken into
account. Indeed, in this paper we propose a dual approach to handle both.
The optimization set-up in \cite{wei20131} is similar to the one considered in 
this paper. Differently from our approach, which is a dual method, a primal-dual
algorithm is proposed in this reference. This difference results in different 
algorithm as well as different requirements on the cost functions.

 The contribution of the paper is twofold. First, for a fixed graph topology, we
 develop a distributed optimization algorithm (based on a centralized dual
 proximal gradient idea introduced in Beck \cite{beck2009fast}) to minimize a
 separable strongly convex cost function. 
 The proposed distributed algorithm is based on a proper choice of primal
 constraints (suitably separating the graph-induced and node-local constraints),
 that gives rise to a dual problem with a separable structure when expressed in
 terms of local conjugate functions. Thus, a proximal gradient applied to such a
 dual problem turns out to be a distributed algorithm where each node updates:
 (i) its primal variable through a local minimization and (ii) its dual
 variables through a suitable local proximal gradient step.  The algorithm
 inherits the convergence properties of the centralized one and thus exhibits an
 $O(1/t)$ rate of convergence in objective value. We point out that the
 algorithm can be easily accelerated through a Nesterov's scheme,
 \cite{nesterov2007gradient}, thus obtaining an $O(1/t^2)$ rate.

Second, we propose an asynchronous version of this algorithm for a symmetric
gossip communication protocol. 
In this \emph{event-triggered} communication 
set-up, a node is in idle mode until its local timer triggers.
When in idle, it collects messages from neighboring nodes that are awake and may
send information if required. When the local timer triggers, it updates its
local (primal and dual) variables and broadcasts them to neighboring
nodes. Under mild assumptions on the local triggering timers, the whole
algorithm results into a random choice of one active node per
iteration. Convergence is proven by showing that the distributed algorithm
corresponds to a block-coordinate proximal gradient, as the one proposed in
\cite{richtarik2014iteration}, performed on the dual problem.
An important feature of the distributed algorithm is that each node can use its
own local step-size, based on the Lipschitz constant of its and its neighbors'
local cost functions. A key distinctive feature of the algorithm analysis is the
combination of duality theory, coordinate-descent methods and properties of the
proximal operator when applied to conjugate functions.

The paper is organized as follows. In Section~\ref{sec:set-up_network} we set-up
the optimization problem and the network model. In
Section~\ref{sec:algorithm_fixed} we introduce and analyze a distributed dual
proximal gradient algorithm for fixed communication graphs, while in
Section~\ref{sec:algorithm_gossip} we extend the algorithm to an asynchronous
scenario. In Section~\ref{sec:simulations} we show through a numerical
example the convergence properties of the asynchronous algorithm.

Due to space constrains all proofs are omitted in this paper and will be
provided in a forthcoming document.
\paragraph*{Notation}
Given a closed, nonempty convex set $X$, the indicator function of $X$ is
defined as $I_X(x) = 0$ if $x\in X$ and $I_X(x) = +\infty$ otherwise.
Let $f \,:\, \real^d \to \real \cup \{+\infty \}$, its conjugate function $f^*\,:\, \real^d \to
\real$ is defined as $f^*(y) := \sup_x \left \{ y^Tx - f(x)\right \}$.
Let $f \,:\, \real^d \to \real \cup \{+\infty \}$ be a closed proper convex
function and $\alpha$ a positive scalar, the proximal operator $\prox_{\alpha f}
\,:\, \real^d \to \real^d$ is defined by $\prox_{\alpha f} (v) := \argmin_x
\left \{f(x) + \frac{1}{2\alpha} \| x-v\|^2 \right \}$.

\section{Problem set-up and network model}
\label{sec:set-up_network}

We consider the following optimization problem
\begin{align*}
  \min_x \sum_{i=1}^n f_i(x) + g_i(x),
\end{align*}
where $f_i : \real^d \to \real \cup \{+\infty \}$ are proper, closed and strongly
convex extended real-valued functions with strong convexity parameter
$\sigma_i>0$ and $g_i : \real^d \to \real\cup \{+\infty \}$ are proper, closed and
convex extended real-valued functions.
The next assumption is standard and will guarantee that the dual problem
is feasible and equivalent to the primal one (strong duality).

\begin{assumption}[Constraint qualification]
  The~intersection of the relative interior of $\dom\,\sum_{i=1}^n f_i$ and the relative interior
  of $\dom\, \sum_{i=1}^n g_i$ is non-empty. \oprocend

\label{ass:Slater}
\end{assumption}

Notice that a convex constrained optimization problem
\begin{align}
	\notag
	\min_x & \; \sum_{i=1}^n f_i ( x) \\
	\notag
	\subj & \; x \in \bigcap_{i=1}^n X_i \subseteq \real^d.
\end{align}
where $X_i$ are convex set, can be modeled in our framework by setting
$g_i(x) = I_{X_i}(x)$.

We want this optimization problem to be solved by a network of processors in a
distributed way, i.e., by a set of peer processors communicating asynchronously
and without the presence of a central coordinator.

Formally, we consider a network of nodes $\until{n}$
communicating according to an asynchronous broadcast protocol.
Each node has its own concept of time defined by a local timer that randomly and
independently of the other nodes triggers when to awake itself. Between two
triggering events the node is in an \emph{idle} mode, i.e., it can receive
messages from neighboring nodes. When a trigger occurs, it switches into an
\emph{awake} mode in which it updates its local variables and transmits the
updated information to its neighbors.

We assume that the asynchronous communication occurs among nodes that are
neighbors in a given fixed, undirected and connected graph $\GG = (\until{n},\EE)$,
where $\EE\subseteq \until{n} \times \until{n}$ is the set of edges. That is, the edge
$(i,j)$ models the fact that node $i$ can receive (respectively send)
information from (to) node $j$ when in idle (awake) mode.  We denote by
$\nbrs_i$ the set of \emph{neighbors} of node $i$ in the fixed graph $\GG$,
i.e., $\nbrs_i := \left\{j \in \until{n} \mid (i,j) \in \EE \right\}$, and by
$|\nbrs_i|$ its cardinality.

We make the following assumption on the local timers.
\begin{assumption}[Exponential i.i.d. local timers]%
  The~waiting times between consecutive triggering events, $T_i$, $i\in
  \until{n}$, are exponential i.i.d.\ random variables.  \oprocend
  \label{ass:timers}
\end{assumption}
Let $i_t \in \until{n}$, $t = 1, 2, \ldots$ be the sequence identifying the
generic $t$-th triggered node. Assumption~\ref{ass:timers} implies that $i_t$ is
an i.i.d.\ uniform process on the alphabet $\until{n}$. Each triggering will
induce an iteration of the distributed optimization algorithm, so that $t$ will
be a universal, discrete time indicating the $t$-th iteration of the algorithm
itself.

To exploit the sparsity of the underlying graph, we
introduce copies of $x$ and a coherence consensus constraint, so that the
optimization problem can be equivalently written as
\begin{align}
	\begin{split}
	\min_{x_1, \ldots, x_n} &\; \sum_{i=1}^n  f_i (x_i) + g_i (x_i) \\
	\subj & \; x_i = x_j \hspace{0.4cm} \forall\, (i,j) \in \EE
	\end{split}
\label{eq:primal_problem_x}
\end{align}
with $x_i\in \real^d$ for all $i\in\until{n}$. The connectedness of
$\mathcal{G}$ guarantees the equivalence.

\section{Algorithm for fixed communication graph}
\label{sec:algorithm_fixed}
In this section we derive and analyze a distributed dual proximal gradient
algorithm for a fixed graph.

\subsection{Dual problem derivation}
We derive the dual version of the problem that will allow us to design our
distributed dual proximal gradient algorithm. To obtain the desired separable structure of
the dual problem, we set-up an equivalent formulation of
problem~\eqref{eq:primal_problem_x} by adding a new set of slack variables $z_i$,
$i\in\until{n}$, i.e., 
\begin{align}
	\begin{split}
	\min_{\stackrel{x_1, \ldots, x_n}{z_1, \ldots, z_n}} & \; \sum_{i=1}^n  f_i (x_i) + g_i (z_i) \\
	\subj & \; x_i = x_j \hspace{0.4cm} \forall\, (i,j) \in \EE \\
	& \; x_i = z_i \hspace{0.5cm} \forall\, i \in \until{n}.
	\end{split}
\label{eq:primal_problem}
\end{align}

Let ${\bf x} = [x_1^T \,\ldots\, x_n^T]^T$ and ${\bf z}  = [z_1^T \, \ldots\, z_n^T]^T$,
the Lagrangian of the primal problem~\eqref{eq:primal_problem} is given by
\begin{align*}
	L({\bf x},{\bf z},\Lambda,\mu) & = \sum_{i=1}^n \bigg( f_i (x_i) + g_i (z_i) 
	\\
	& \hspace{1cm} + \sum_{j\in \nbrs_i} \Big(\lambda_i^j \Big)^T
	(x_i-x_j) + \mu_i^T(x_i  - z_i) \bigg)
	\\
  & = \sum_{i=1}^n \bigg( f_i (x_i) + \sum_{j\in \nbrs_i} \Big( \lambda_i^j \Big)^T (x_i-x_j) + \mu_i^Tx_i
  \\
  & \hspace{1cm} + g_i (z_i) - \mu_i^Tz_i \bigg),
\end{align*}
where $\Lambda$ and $\mu$ are respectively the vectors of the Lagrange
multipliers $\lambda_i^j$, $(i,j)\in\EE$, and $\mu_i$, $i\in\until{n}$, and in
the last line we have separated the terms in $\textbf{x}$ and $\textbf{z}$.
Since $\GG$ is undirected, the Lagrangian can be equivalently
rewritten as 
\begin{align*}
  L({\bf x},{\bf z},\Lambda,\mu) & =\sum_{i=1}^n \bigg( f_i (x_i) + x_i^T 
      \bigg ( \sum_{j\in \nbrs_i} (\lambda_i^j -\lambda_j^i) + \mu_i \bigg ) 
  \\
  & \hspace{1cm} + g_i (z_i) - z_i^T\mu_i \bigg)
\end{align*}
where $\lambda_i^j$, $j \in \nbrs_i$ and $\mu_i$ are variables handled by node
$i$ (consistently $\lambda_j^i$ is handled by node $j$ neighbor of node $i$).

The dual function is
\begin{align*}
  q(\Lambda,\mu) & := \min_{{\bf x,z}} L({\bf x,z},\Lambda, \mu) 
  \\
  & = \min_{\bf x} \sum_{i=1}^n \bigg( f_i (x_i) + x_i^T 
      \bigg( \sum_{j\in \nbrs_i} (\lambda_i^j -\lambda_j^i) + \mu_i \bigg) \bigg) 
  \\
  & \hspace{1cm} + \min_{\bf z} \sum_{i=1}^n \Big(g_i (z_i) - z_i^T\mu_i \Big)
  \\
  & = \sum_{i=1}^n \min_{x_i} \bigg( f_i (x_i) + x_i^T 
      \bigg( \sum_{j\in \nbrs_i} (\lambda_i^j -\lambda_j^i) + \mu_i \bigg) \bigg) 
  \\
  & \hspace{1cm} + \sum_{i=1}^n \min_{z_i} \Big(g_i (z_i) - z_i^T\mu_i \Big)
\end{align*}
where we have used the separability of the Lagrangian with respect to each $x_i$
and each $z_i$. Then, by using the definition of conjugate function, the dual
function can be rewritten as
\begin{align*}
  q(\Lambda,\mu) & = \sum_{i=1}^n \bigg( -f_i^* \bigg(-\sum_{j\in \nbrs_i}
      (\lambda_i^j -\lambda_j^i) - \mu_i \bigg) - g_i^*(\mu_i)\bigg).
\end{align*}

The dual problem of \eqref{eq:primal_problem} consists of maximizing the dual function 
with respect to dual variables $\Lambda$ and $\mu$, i.e.,
\begin{align}
\label{eq:dual_problem}
	\!\!\max_{\Lambda,\mu} &\; 
	\sum_{i=1}^n \bigg( -f_i^* \bigg (-\sum_{j\in \nbrs_i} (\lambda_i^j -\lambda_j^i) - \mu_i \bigg ) - g_i^*(\mu_i) \bigg).
\end{align}

By Assumption~\ref{ass:Slater} the dual problem~\eqref{eq:dual_problem} is
feasible and strong duality holds, so that \eqref{eq:dual_problem} can be solved
to get a solution of \eqref{eq:primal_problem}.

\subsection{Distributed Dual Proximal Gradient Algorithm}
To develop the algorithm, we start rewriting problem~\eqref{eq:dual_problem} by
using a more compact notation, and in the equivalent minimization version.
We stack the dual variables as $y = [y_1^T \;\ldots\; y_n^T]^T$, where
\begin{align}
  \label{eq:yi}
  y_i = 
	\begin{bmatrix}
		\Lambda_i\\ 
		\mu_i
  \end{bmatrix}
	\in \real^{d |\nbrs_i|+d}
\end{align}
with $\Lambda_i\in \real^{d|\nbrs_i|}$ a vector whose block-component associated to
neighbor $j$ is $\lambda_i^j \in \real^d$.
Thus, the dual problem can be written as
\begin{align}
   \label{eq:dual_min_problem}
	\min_{y} &\; \Gamma(y) = F^*(y) + G^*(y),
\end{align}
where 
\begin{align}
  \notag
  F^*(y) & := \sum_{i=1}^n  f_i^* \Big(-\sum_{j\in \nbrs_i} (\lambda_i^j -\lambda_j^i) - \mu_i \Big) \\
  \notag
  G^*(y) & := \sum_{i=1}^n g_i^* \big( \mu_i \big).
\end{align}

The proposed distributed algorithm will be based on a proximal gradient applied
to the above formulation of the dual problem. Next, we describe the local update
of each node $i\in\until{n}$ and then, in the next subsection, we show its
convergence properties.

Node $i$ updates its local dual variables $\lambda_i^j$, $j \in \nbrs_i$, and
$\mu_i$ according to a local proximal gradient step, and its primal
variable $x_i^\star$ through a local minimization. The step-size of the proximal
gradient step is denoted by $\alpha$. Then, the updated primal and dual
values are exchanged with the neighboring nodes according to a synchronous
communication over a fixed undirected graph. 
The local dual variables at node $i$ are initialized as $\lambda_{i0}^j$,
$j\in\nbrs_i$, and $\mu_{i0}$.
A pseudo-code of the local update at each node of the distributed algorithm is
given in Algorithm~\ref{alg:fixed}.%
\begin{algorithm}
  \begin{algorithmic}[0]
    \State {\bf Processor states:} $x_i^\star$, $\lambda_i^j$ for all $j\in \nbrs_i$
    and $\mu_i$
    \State {\bf Initialization:}
    $\lambda_i^j(0)=\lambda_{i0}^j$
    for all $j\in \nbrs_i$, $\mu_i(0)=\mu_{i0}$ 
    \StatexIndent[1] 
    \begin{small}
    $x_i^\star(0) =\argmin_{x_i} \left \{ \!x_i^T\! \left (\! \sum_{j\in \nbrs_i} \! \left (\!
             \lambda_{i0}^j - \lambda_{j0}^i \!\right ) \!+\! \mu_{i0} \!\right )
         \! + \! f_i(x_i)  \!\right \}$
    \end{small} 

    \State {\bf Evolution:}
    \StatexIndent[0.3] \textsc{for:} $t=1,2,\ldots$

    \StatexIndent[0.6] receive $x_j^\star(t-1)$ for each $j\in \nbrs_i$
    \StatexIndent[0.6] update
      \begin{align}
        \notag \!\!\!\lambda_i^j(t) & = \lambda_i^j(t-1) + \alpha \big [ x_i^\star(t-1) -
          x_j^\star(t-1) \big]
      \end{align}
    \StatexIndent[0.6] update
      \begin{align*}
        & \tilde{\mu}_i = \mu_i(t-1) + \alpha \; x_i^\star(t-1) 
        \\
        & \mu_i(t) = \prox_{\alpha g_i^*}  \big( \tilde{\mu}_i \big) =
        \tilde{\mu}_i - \alpha\;\prox_{\frac{1}{\alpha}g_i} \bigg(
        \frac{\tilde{\mu}_i}{\alpha} \bigg)
	  \end{align*}
    \StatexIndent[0.6] receive $\lambda_j^i(t)$ for each $j\in \nbrs_i$
    \StatexIndent[0.6] update
      \begin{small}
      \begin{align}
      \notag
        \!\!\! x_i^\star(t)
       = \argmin_{x_i} \bigg\{  x_i^T \bigg( \sum_{j\in \nbrs_i} \! 
       \Big (\! \lambda_i^j(t) - \lambda_j^i(t) \! \Big ) \!+\! \mu_i(t) \! \bigg)
         \! + \! f_i(x_i) \bigg\}
      \end{align}      
      \end{small}
  \end{algorithmic}

\caption{\small \DDPGfixed} \label{alg:fixed}
\end{algorithm}
\begin{remark}
In order to start the algorithm, a preliminary communication step is needed in
which each node $i$ sends to each neighbor $j$ its $\lambda_{i0}^j$. This step
can be avoided if the nodes agree to set $\lambda_{i0}^j =0$.~\oprocend
\end{remark}

\subsection{Algorithm analysis}

\begin{lemma}[\cite{boyd2004convex,beck2014fast}]
  Let $\varphi$ be a closed, strictly convex function and $\varphi^*$ its conjugate function.
  Then 
  \begin{align}
    \notag
    \nabla \varphi^*(y) \!=\! \argmax_x \left\{ y^Tx \!-\!\varphi(x)\right\} 
    \!=\!  \argmin_x \left\{ \varphi(x) \!-\! y^Tx \right\}\!. 
  \end{align}
  Moreover, if $\varphi$ is strongly convex with convexity parameter $\sigma$, then 
  $\nabla \varphi^*$ is Lipschitz continuous with Lipschitz constant given by 
  $\frac{1}{\sigma}$. \oprocend  
\label{lem:conjugate_gradient}
\end{lemma}

\begin{lemma}[Moreau decomposition, \cite{parikh2013proximal}]
	Let $f : \real^d \rightarrow \real\cup \!\{+\infty \}$ be a closed,
  strictly convex function and $f^*$ its conjugate, then 
    $\forall x\in \real^d$, $x = \prox_f(x) + \prox_{f^*}(x)$. \oprocend
	\label{lem:Moreau_decomposition}
\end{lemma}

\begin{lemma}[Extended Moreau decomposition]
\label{lem:extended_Moreau}
Let $\varphi : \real^d \rightarrow \real\cup \{+\infty \}$ be a closed,
strictly convex function and $\varphi^*$ its conjugate. Then for any $x\in
\real^d$ and $\alpha>0$, it holds $x = \prox_{\alpha \varphi}\left\{ x\right\} +
\alpha \prox_{\frac{1}{\alpha} \varphi^*} \left\{ \frac{x}{\alpha} \right\}$.\oprocend
\end{lemma}

\begin{lemma}
  Let $y = [y_1^T \; \ldots \; y_n^T]^T \in \real^{n(D+d)}$ where
  $y_i = [\Lambda_i^T\;\mu_i^T]^T$ with $\Lambda_i \in \real^D$ and $\mu_i \in\real^d$,
  $i\in\until{n}$.  Let $G^*(y) = \sum_{i=1}^n g_i^*(\mu_i)$, then the proximal operator of $\alpha G^*$ 
  evaluated at $y$ is given by
  \begin{align}
    \notag
    \prox_{\alpha G^*} \big( y \big)
    & =
    \begin{bmatrix}
      \Lambda_1 \\
      \prox_{\alpha g_1^*}( \mu_1)
      \\ \vdots \\
      \Lambda_n \\
      \prox_{\alpha g_n^*}( \mu_n)
    \end{bmatrix}.
  \end{align}%
  \oprocend
  \label{lem:proxG}
\end{lemma}

\begin{theorem}
  For each $i\in\until{n}$, let $f_i$ be a proper, closed, strongly convex
  extended real-valued function with strong convexity parameter $\sigma_i>0$,
  and let $g_i$ be a proper, convex extended real-valued function. Let $y^\star$ be
  the minimizer of \eqref{eq:dual_min_problem}.  Suppose that in
  Algorithm~\ref{alg:fixed} the step-size $\alpha$ is chosen such that $0<\alpha
  \leq \dfrac{1}{\sum_{i=1}^n \frac{1}{\sigma_i}}$.  Then the sequence $y(t) =
  [y_1(t)^T \ldots y_n(t)^T]^T$ generated by the \DDPGfixed\,
  (Algorithm~\ref{alg:fixed}) converges to $y^\star$ and in objective value
  satisfies
	\begin{align}
    	\notag
	  \Gamma(y(t)) - \Gamma(y^\star) \leq \frac{\Big(\sum_{i=1}^n \frac{1}{\sigma_i}\Big) \| y_0 - y^\star\|^2}{2t},
	\end{align}
  where $y_0 = [y_1(0)^T \; \ldots \; y_n(0)^T]^T$ is the initial condition.

\label{thm:fixed}
\end{theorem}

\begin{remark}[Nesterov's acceleration]
  We can include a Nesterov's \emph{extrapolation step} in the algorithm, which
  accelerates the algorithm (\cite{nesterov2007gradient} for further details),
  attaining a faster $O(1/t^2)$ convergence rate in objective value.
  \oprocend
\end{remark}

\section{Asynchronous distributed dual proximal gradient}
\label{sec:algorithm_gossip}
In this section we present an asynchronous distributed dual proximal gradient
and prove its convergence in probability.

We start by describing the local evolution at each node $i\in\until{n}$. 
First, recall from the network model introduced in
Section~\ref{sec:set-up_network} that a node can be into two different modes:
when in \emph{idle} it continuously listens to incoming messages from its
neighbors (and, if needed, may send them auxiliary information back),
while when in \emph{awake} it updates its local variables and transmits them 
to its neighbors.
The transition between modes is asynchronously ruled via local timers,
$\tau_i\in\real$, $i\in \until{n}$ (they are assumed to have infinite
precision). As from Assumption~\ref{ass:timers}, timers trigger according
to $n$ exponential i.i.d. random variables $T_i$, $i\in \until{n}$. In the
algorithm we make a slight abuse of notation denoting by $T_i$ the realization
of the random variables $T_i$.

Each node $i$ updates its local dual variables $\lambda_i^j$, $j\in\nbrs_i$ and
$\mu_i$ by a local proximal gradient step, and its primal variable
$x_i^\star$ through a local minimization. Each node uses a properly chosen,
\emph{local} step-size $\alpha_i$ for the proximal gradient step.

\begin{algorithm}
  \begin{algorithmic}[0] 
    \Statex {\bf Processor states:} 
    $x_i^\star$, $\lambda_i^j$ for all $j\in \nbrs_i$ and $\mu_i$

    \State {\bf Initialization:} $\lambda_i^j=\lambda_{i0}^j$ for all $j\in
    \nbrs_i$, $\mu_i=\mu_{i0}$ and %
    \StatexIndent[1] 
    \begin{small}
      $\!\!x_i^\star(0) =\argmin_{x_i} \left \{ \!x_i^T\! \left (\! \sum_{j\in \nbrs_i} \! \left (\!
             \lambda_{i0}^j - \lambda_{j0}^i \!\right ) \!+\! \mu_{i0} \!\right )
         \! + \! f_i(x_i)  \!\right \}$
    \end{small}
    \StatexIndent[1] set $\tau_i = 0$ and get $T_i$ \medskip
    
    \Statex {\bf Evolution:} 
    \Label \texttt{\textbf{\textit{IDLE:}}}

    \StatexIndent[0.5]\textsc{while}{ $\tau_i \leq T_i$} %
    \StatexIndent[1] receive $x_j^\star$ and/or $\lambda_j^i$ from each
    $j\in \nbrs_i$  %
    \StatexIndent[1] \textsc{if} $\lambda_j^i$ is received \textsc{then} compute and broadcast
      \begin{small}
      \begin{align}
      \notag
      {x_i^\star}
        = \argmin_{x_i} \bigg\{ x_i^T \bigg( \sum_{\ell\in \nbrs_{i}} \bigg(
        \lambda_{i}^\ell - \lambda_\ell^{i} \bigg)+ \mu_{i} \bigg)
        + f_{i}(x_i) \bigg\}
      \end{align}
      \end{small}

         \Label \texttt{\textbf{\textit{AWAKE:}}}
  
   \StatexIndent[0.5] update and broadcast
      \begin{align*}
        \hspace{-1.1cm}
        {\lambda_{i}^j}^+ & = \lambda_{i}^j + \alpha_i \big ( x_i^\star - x_j^\star \big ),\:\: \forall \, j \in \nbrs_{i}
        \end{align*}
   \StatexIndent[0.5] update
    \begin{align*}
      & \tilde{\mu}_{i}  = \mu_{i} + \alpha_i \; x_{i}^\star 
      \\
      & \mu_{i}^+  = \prox_{\alpha_i g_{i}^*} \big( \tilde{\mu}_{i} \big)
      = \tilde{\mu}_{i} - \alpha_i\;\prox_{\frac{1}{\alpha_i} g_{i} } \Big(\frac{\tilde{\mu}_{i}}{\alpha_i} \Big)
    \end{align*}

   \StatexIndent[0.5] compute and broadcast
      \begin{small}
      \begin{align}
      \notag
      {x_i^\star}
        = \argmin_{x_i} \bigg\{ x_i^T \bigg( \sum_{j\in \nbrs_{i}} \bigg(
        {\lambda_{i}^j}^+ - \lambda_j^{i} \bigg)+ \mu^+_{i} \bigg)
        + f_{i}(x_i) \bigg\}
      \end{align}
      \end{small}

   \StatexIndent[0.5] set $\tau_i=0$, get a new $T_i$ and go to \texttt{\textbf{\textit{IDLE}}}.

  \end{algorithmic}
  \caption{\small \DDPGgossip}
  \label{alg:gossip}
\end{algorithm}
\begin{remark}
  In order to set the step-size $\alpha_i$, node $i$ needs a preliminary
  communication step to receive the convexity parameters from its
  neighbors.~\oprocend
\end{remark}
From an external, global perspective, the described local asynchronous updates
result into an algorithmic evolution, in which at each iteration only one node
wakes up randomly, uniformly and independently from previous iterations. This
follows from the memoryless property of the exponential distribution. Thus, in
this high-level view, we can consider a \emph{universal} (discrete)
time-variable $t$, which counts the iterations of the whole algorithm
evolution. 
This variable will be used in the statement 
of Theorem~\ref{thm:gossip}.

\begin{theorem}
  For each $i\in\until{n}$, let $f_i$ be a proper, closed and
  strongly convex extended real-valued function with strong convexity parameter
  $\sigma_i>0$, and let $g_i$ be a proper convex extended real-valued
  function. Let $y^\star$ be the minimizer of \eqref{eq:dual_min_problem}.
  Suppose that in Algorithm~\ref{alg:gossip} each local step-size $\alpha_i$ 
  is chosen such that $0<\alpha_i \leq \frac{1}{L_i}$ with
  \begin{align}
    \notag
    L_i = \sqrt{ \frac{1}{\sigma_i^2}+ \sum_{j\in \nbrs_i}
        \Big(\frac{1}{\sigma_i} + \frac{1}{\sigma_j}\Big)^2 }.
  \end{align}
  Then the sequence $y(t) = [y_1(t)^T \ldots y_n(t)^T]^T$ generated by the
  \DDPGgossip\, (Algorithm~\ref{alg:gossip}) converges in probability to $y^\star$,
  i.e., for any $\varepsilon \in \big( 0, \Gamma(y_0) \big)$, where $y_0 = [y_1(0)^T \ldots
  y_n(0)^T]^T$ is the initial condition, and target confidence $0<\rho<1$, there exists
  $\bar{t}(\varepsilon,\rho)>0$ such that for all $t\geq \bar{t}$ it holds
  \begin{align}
    \notag
    \mathrm{P} \Big( \Gamma(y(t)) - \Gamma(y^\star) \leq \varepsilon \Big) \geq 1- \rho.
  \end{align}

\label{thm:gossip}
\end{theorem}

\section{Simulations}
\label{sec:simulations}
In this section we provide a numerical example showing the effectiveness of
the proposed \DDPGgossip.

We consider an undirected connected Erd\H{o}s-R\'enyi
graph $\GG$ with parameter $0.2$, connecting $n=15$ nodes.
We assume each decision variable $x_i \in \real^{2}$, $i \in \until{n}$. 
Let each local objective function $f_i$ be quadratic and randomly generated as
\begin{align}
  \notag
  f_i(x_i) = x_i^TQ_ix_i + r_i^Tx_i
\end{align}
where $Q_i \in \real^{2\times 2}$ is diagonal with diagonal elements uniformly distributed
in $[1,2]$ and $r_i \in \real^{2}$ has elements uniformly randomly distributed in $[-5,5]$.
We let each $g_i$ be the indicator function of a convex polytope $X_i
= \{ x_i \in \real^2\mid a_i^T x_i \leq b_i \}$, with components of $a_i$
generated uniformly in $[\,0,10]$ and components of $b_i$ in $[-5,5]$.
We initialize to zero the dual variables $\lambda_i^j$, $j \in \nbrs_i$, and
$\mu_i$ for all $i \in \until{n}$, and use a constant step-size
$\alpha_i=1$ for all nodes.

\begin{figure}[!htbp]
\centering
  \begin{minipage}[t]{0.22\textwidth}
    \centering
    \includegraphics[scale=0.35]{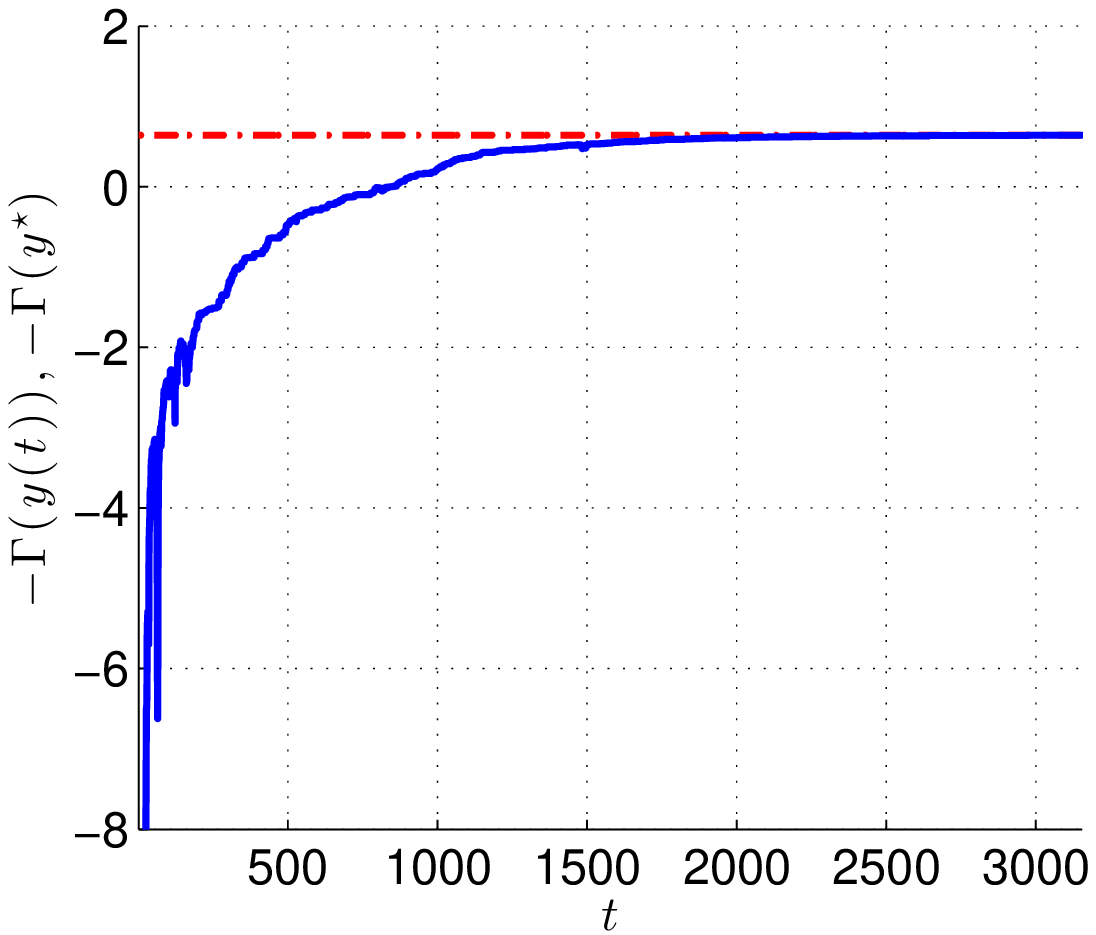}
    \caption{Cost $-\Gamma(y(t))$ (solid blue) vs optimal cost $-\Gamma(y^\star)$ (dotted red).}
    \label{fig:cost}
  \end{minipage}\quad
  \begin{minipage}[t]{0.23\textwidth}
    \centering
    \includegraphics[scale=0.35]{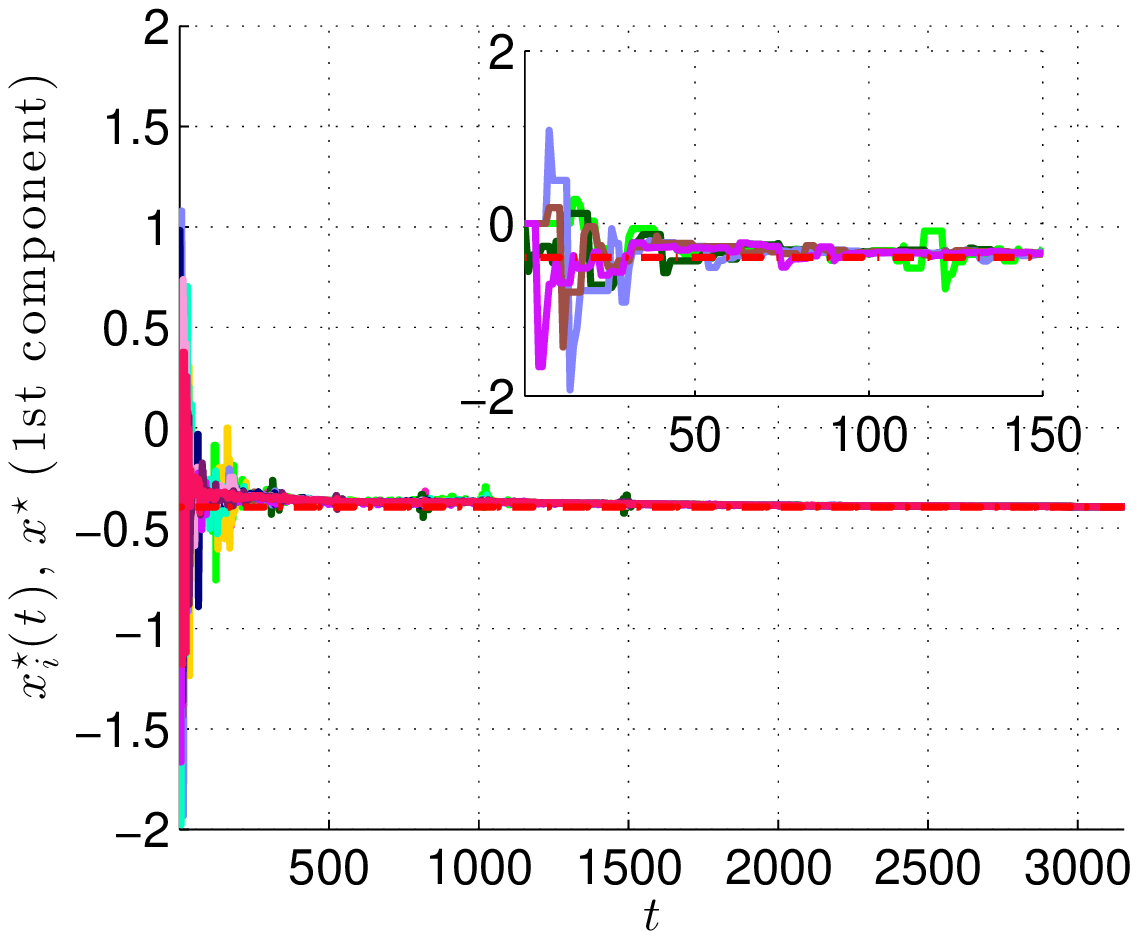}
    \caption{First component of $x_i^\star(t)$, $i \!\in\! \until{n}$, zoom on selected $i$.}
    \label{fig:xi_first_zoom}
  \end{minipage}
\end{figure}

Figure~\ref{fig:cost} shows the convergence of the primal (and dual) cost to the
optimal centralized value. We recall that the primal cost is $-\Gamma(y(t))$,
with $\Gamma(y(t))$ being the dual cost in the minimization version
\eqref{eq:dual_min_problem}.
In Figure~\ref{fig:xi_first_zoom} 
we plot the behavior of the first component of primal variables $x_i^\star(t)$. 
The horizontal dotted-line is the optimal primal solution.
In the inset the first iterations for five selected nodes, $x_i^\star$,
$i=1,5,6,7,13$, are highlighted, in order to better show the transient,
piece-wise constant behavior due to the gossip update.

\begin{figure}[h]
\centering
  \begin{minipage}[t]{0.22\textwidth}
    \centering
	  \includegraphics[scale=0.35]{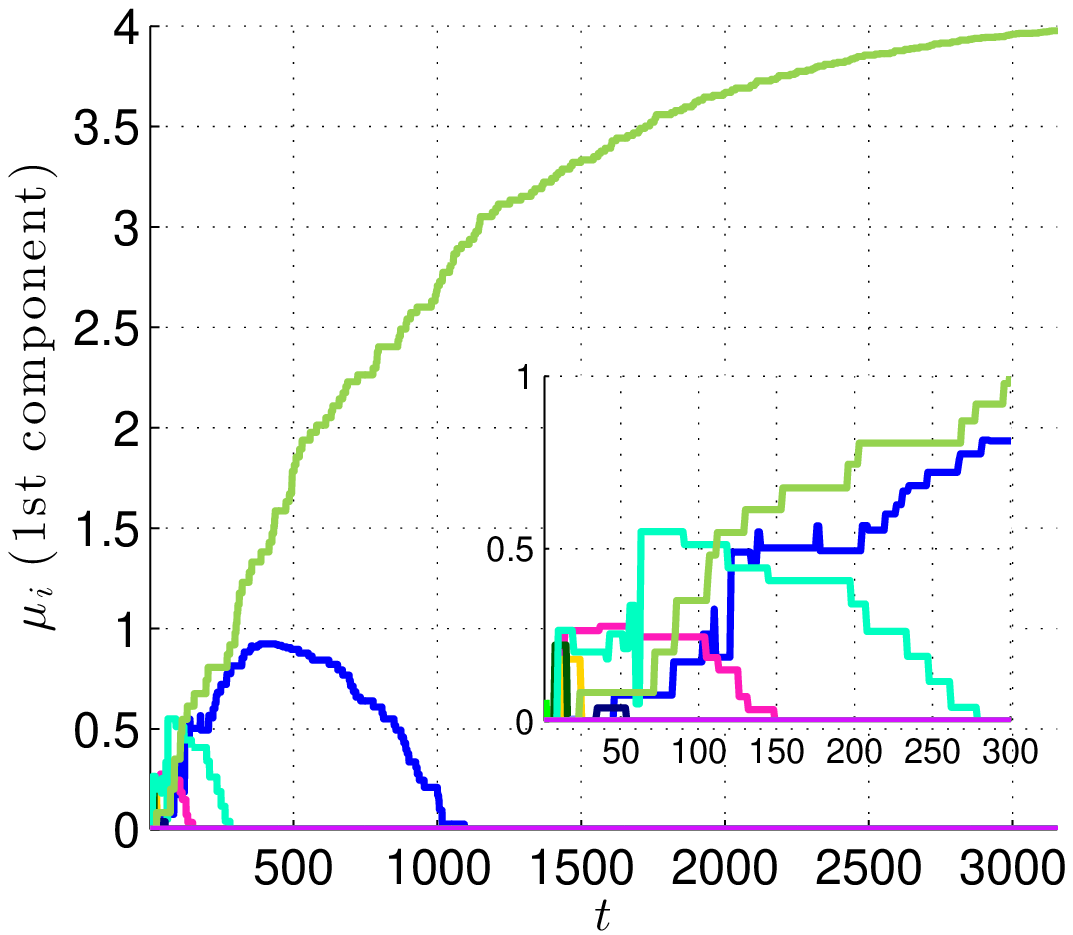}
    \caption{First component of $\mu_i(t)$, $i\in\until{n}$.}
    \label{fig:mu_first}
  \end{minipage}\quad
  \begin{minipage}[t]{0.23\textwidth}
    \centering
    \includegraphics[scale=0.35]{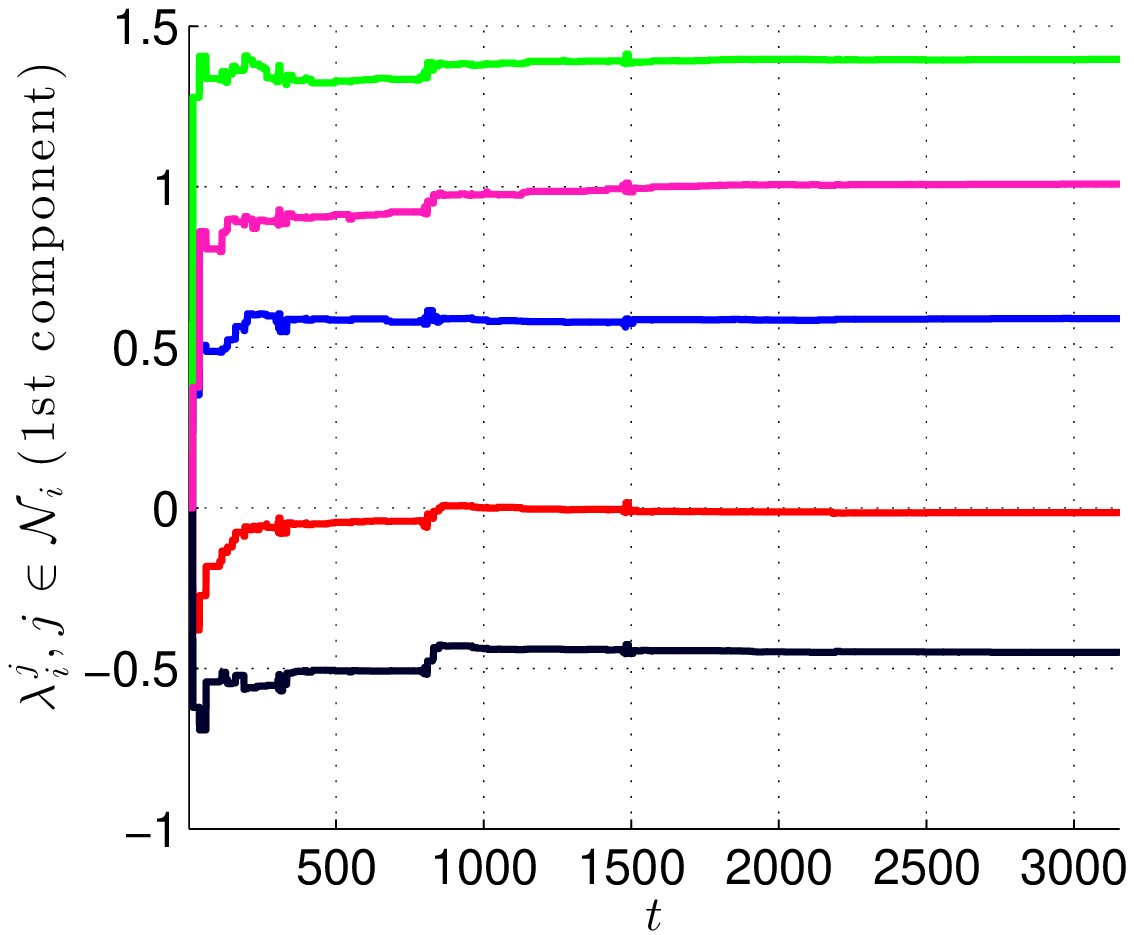}
    \caption{First component of $\lambda_i^j$, $j\in N_i$, with $i=5$.}
    \label{fig:lambda_first}
  \end{minipage}
\end{figure}

Then we show the evolution of the dual variables.
First, note that $\mu_i$ is associated to the local constraint $X_i$ of
$x_i$. We obtain that only $\mu_{13}$, the multiplier relative to the only active
constraint, converges to a nonzero value, whereas all the other $\mu_i$s,
associated to the inactive constraints, converge to $0$. 
In Figure~\ref{fig:mu_first} the first component of $\mu_{13}$ is plotted.

Finally, we plot the evolution of $\lambda_i^j$, $j\in\nbrs_i$, for node $i=5$
(with $\nbrs_i = \{3,6,10,12,14\}$), see Figure~\ref{fig:lambda_first} for the
first component. As expected the multipliers converge to nonzero values
representing the ``price'' needed to enforce equality constraints on the primal
variables $x_i$ and $x_j$.

\section{Conclusions}
\label{sec:Conclusions}
In this paper we have proposed an asynchronous, distributed optimization
algorithm, based on a block-coordinate dual proximal gradient method to solve
separable, constrained optimization problems.
The main idea is to construct a suitable, separable dual problem via a proper
choice of primal constraints. Then, the dual problem is solved through a
proximal gradient algorithm.
Thanks to the separable structure of the dual problem in terms of local
conjugate functions, the proximal gradient update results into a distributed
algorithm, where each node performs a local minimization on its primal variable,
and a local proximal gradient update on its dual variables.
An asynchronous version of the distributed algorithm is obtained by exploiting a
randomized, block-coordinate descent approach.

\appendix

\renewcommand{\thetheorem}{\thesection.\arabic{theorem}}
\renewcommand{\theequation}{\thesection.\arabic{equation}}

\section{Randomized coordinate descent for composite functions}

Consider the following optimization problem
\begin{align}
  \label{eq:appendix_prob}
	\min_{y\in \real^N}\; \Gamma(y) :=	\Phi(y) + \Psi(y)
\end{align}
where $\Phi : \real^N \to \real$ and $\Psi : \real^N \to \real \cup
\{+\infty\}$ are convex functions.

We decompose the decision variable as $y = [y_1^T \; \ldots \; y_n^T]^T$ and,
consistently, we decompose the space $\real^N$ into $n$ subspaces as follows.
Let $U \in \real^{N\times N}$ be a column permutation of the $N \times N$
identity matrix and, further, let $U = [U_1 \; U_2\;\ldots\;U_n]$ be a
decomposition of $U$ into $n$ submatrices, with $U_i\in \real^{N\times N_i}$ and
$\sum_i N_i = N$.  Thus, any vector $y\in \real^N$ can be uniquely written as
$y = \sum_i U_iy_i$ and, viceversa, $y_i = U_i^Ty$.

We let problem~\eqref{eq:appendix_prob} satisfy the following assumptions.

\begin{assumption}[\textbf{Smoothness of $\Phi$}]%
  The gradient of $\Phi$ is block coordinate-wise
  Lipschitz continuous with positive constants $L_1,\ldots,L_n$. That is, for
  all $y \in \real^N$ and $s_i\in \real^{N_i}$ it holds
  \begin{align*}
    \|\nabla_i \Phi(y + U_i s_i) - \nabla_i \Phi(y) \| \leq L_i \|s_i\|,
  \end{align*}
  where $\nabla_i \Phi(y)$ is the $i$-th block component of $\nabla \Phi(y)$.~\oprocend

\label{ass:smoothness}
\end{assumption}
\begin{assumption} [\textbf{Separability of $\Psi$}]
  The function $\Psi$ is block-separable, 
  i.e., it can be decomposed as $\Psi(y) = \sum_{i=1}^n\psi_i (y_i)$,
  with each $\psi_i : \real^{N_i} \to \real \!\cup \! \{+\infty\}$ a
  proper, closed convex extended real-valued function.~\oprocend
\label{ass:separability}
\end{assumption}

\begin{assumption}[\textbf{Feasibility}]
  The set of minimizers of problem~\eqref{eq:appendix_prob} is
  non-empty.~\oprocend
\label{ass:feasibility}
\end{assumption}

\begin{algorithm}
  \begin{algorithmic}[0]
  
  \State \textbf{Initialization}: $y(0) = y_0$

  \For{$t=0,1,2,\ldots$}{}
    \StatexIndent[0.5] choose $i_t\in \{1,\ldots,n\}$ with probability $\frac{1}{n}$

    \StatexIndent[0.5] compute
    \begin{align}
    \notag
      T^{(i_t)} \big( y(t) \big) = \argmin_{w_{i_t} \in \real^{N_{i_t}}} \Big\{ V_{i_t}( y(t),w_{i_t} ) \Big\}
    \end{align}    
    \StatexIndent[0.5] where
    \begin{align}
    \notag
      V_{i_t} (y, s_{i_t}) := \nabla_{i_t} \Phi (y)^T s_{i_t} \!+\!\frac{L_{i_t}}{2} \| s_{i_t} \|^2 \!+\! \psi_{i_t}( y_{i_t} \!+\! s_{i_t} )
    \end{align}
    
    \StatexIndent[0.5] update $y(t+1) = y(t) + U_{i_t} T^{(i_t)} \big( y(t) \big)$

  \EndFor
  \end{algorithmic}
  
  \caption{UCDC} \label{alg:ucdc}
\end{algorithm}

The convergence result for UCDC (Algorithm~\ref{alg:ucdc}) is given in \cite[Theorem~5]{richtarik2014iteration},
here reported for completeness.

\begin{theorem}[Theorem~5, \cite{richtarik2014iteration}]
  Let Assumptions~\ref{ass:smoothness}, \ref{ass:separability} and
  \ref{ass:feasibility} hold.  Then, for any
  $\varepsilon \in \Big( 0, \Gamma(y_0) - \Gamma(y^\star) \Big)$, there exists
  $\bar{t}(\varepsilon,\rho)>0$ such that if $y(t)$ is the random sequence
  generated by UCDC applied to problem
  \eqref{eq:appendix_prob}, then for all $t\geq \bar{t}$ it holds that
  \begin{align*}
    \mathrm{P} \Big (\Gamma(y(t)) - \Gamma(y^\star) \leq \varepsilon \Big ) \geq 1- \rho,
  \end{align*}
  where $y^\star$ is a minimizer of problem \eqref{eq:appendix_prob}, 
  $y_0 \in \real^N$ is the initial condition and $\rho \in (0,1)$ is the 
  target confidence.~\oprocend

\label{thm:ucdc}
\end{theorem}

{\small \bibliographystyle{IEEEtran} \bibliography{proximal} }

\end{document}